\documentclass[reqno,10pt,a4paper]{amsart}
\usepackage[margin=2.7cm]{geometry}
\usepackage{cite}
\usepackage{newpxtext}
\usepackage{cabin}
\usepackage[varqu,varl]{inconsolata}
\usepackage[bigdelims,vvarbb]{newpxmath}
\usepackage[cal=cm,bb=esstix,bbscaled=1.05]{mathalfa}
\usepackage{mathtools}
\usepackage{bm}
\usepackage{hyperref}
\usepackage{tikz}
\usetikzlibrary{decorations,decorations.pathmorphing,calc}
\usepackage{subcaption}

\hypersetup{
colorlinks,
linkcolor={red!50!black},
citecolor={blue!50!black},
urlcolor={blue!80!black}
}

\makeatletter
\renewcommand{\p@subfigure}{}
\makeatother
\keywords{Lattice polytope; smooth polytope; Ehrhart theory; triangulation; integer decomposition property.}
\newcommand{\orcid}[1]{\,\resizebox{8px}{!}{\href{https://orcid.org/#1}{\includegraphics{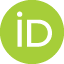}}}}
\DeclareMathOperator{\cone}{cone}
\DeclareMathOperator{\conv}{conv}
\DeclareMathOperator{\vol}{vol}
\newcommand{\Z}{\mathbb{Z}}
\newcommand{\R}{\mathbb{R}}
\newcommand{\orig}{\bm{0}}
\newcommand{\bv}{\bm{v}}
\newcommand{\bw}{\bm{w}}
\newcommand{\bx}{\bm{x}}
\newcommand{\ba}{\bm{a}}
\newcommand{\bb}{\bm{b}}
\newcommand{\be}{\bm{e}}
\newcommand{\cP}{\mathcal{P}}
\newcommand\blfootnote[1]{%
  \begingroup
  \renewcommand\thefootnote{}\footnote{#1}%
  \addtocounter{footnote}{-1}%
  \endgroup
}
\theoremstyle{definition}
\newtheorem{theorem}{Theorem}
\newtheorem{example}{Example}
\newtheorem{question}{Question}
\begin{document}
\author[J.\ Hofscheier]{Johannes Hofscheier\orcid{0000-0001-8642-3984}}
\address{School of Mathematical Sciences, University of Nottingham, Nottingham, NG7 2RD, United Kingdom}
\email{johannes.hofscheier@nottingham.ac.uk}
\author[A.\ M.\ Kasprzyk]{Alexander M.~Kasprzyk\orcid{0000-0003-2340-5257}}
\address{Mathematics Institute, University of Warwick, Coventry, CV4 7AL, United Kingdom}
\email{alexander.kasprzyk@warwick.ac.uk}
\title{Is there a smooth lattice polytope which does not have the integer decomposition property?}
\maketitle
\begin{abstract}
We introduce Tadao Oda's famous question on lattice polytopes which was originally posed at Oberwolfach in 1997 and, although simple to state, has remained unanswered. The question is motivated by a discussion of the two\nobreakdash-dimensional case -- including a proof of Pick's Theorem, which elegantly relates the area of a lattice polygon to the number of lattice points it contains in its interior and on its boundary.
\end{abstract}
\blfootnote{Snapshots of Modern Mathematics from Oberwolfach, SNAP-2025-008-EN (2025), doi:\href{https://doi.org/10.14760/SNAP-2025-008-EN}{10.14760/SNAP-2025-008-EN}}
\section{Introduction}
Lattice polytopes are fundamental objects in mathematics and play a crucial role in a broad range of subjects such as discrete and algebraic geometry, algebra, combinatorics, coding theory, and optimisation theory. They arise naturally in a variety of unexpected or even surprising ways. Consider the following classical question from enumerative combinatorics, for instance.
\begin{question}\label{quest:intro}
How many monomials in three variables of a given degree~$m$ are there?
\end{question}
A~\emph{monomial} is a product of (powers of) certain variables; its~\emph{degree} is the sum of the exponents as they appear in the product. We thus wish to find the total number of all triples of non-negative integers whose sum equals a given integer~$m$.

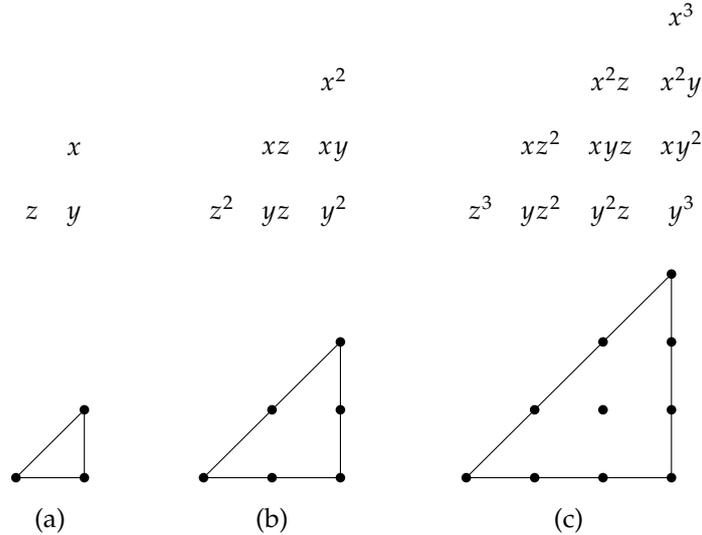
\begin{figure}[t]
\centering
\begin{tabular}{c@{\qquad\qquad}c@{\qquad\qquad}c}
\parbox[b][][t]{1cm}{$\begin{array}{cc}
&x \\\\
z&y 
\end{array}$}&
\parbox[b][][t]{2cm}{$\begin{array}{ccc}
& & x^2 \\\\
&xz & xy \\\\
z^2 & yz & y^2
\end{array}$}&
\parbox[b][][t]{3cm}{$\begin{array}{cccc}
&&&x^3 \\\\
&&x^2z&x^2y\\\\
&xz^2&xyz&xy^2\\\\
z^3&yz^2&y^2z&y^3
\end{array}$}\\\\
\begin{tikzpicture}[scale=0.9]
\draw (0,0) -- (1,0) -- (1,1) -- cycle;
\foreach \x in {0,...,1}
\foreach \y in {0,...,\x}
\fill (\x,\y) circle (2pt);
\end{tikzpicture}&
\begin{tikzpicture}[scale=0.9]
\draw (0,0) -- (2,0) -- (2,2) -- cycle;
\foreach \x in {0,...,2}
\foreach \y in {0,...,\x}
\fill (\x,\y) circle (2pt);
\end{tikzpicture}&
\begin{tikzpicture}[scale=0.9]
\draw (0,0) -- (3,0) -- (3,3) -- cycle;
\foreach \x in {0,...,3}
\foreach \y in {0,...,\x}
\fill (\x,\y) circle (2pt);
\end{tikzpicture}\\[.5em]
(a)&(b)&(c)
\end{tabular}
\caption{Interpreting monomials in three variables~$x,y$, and~$z$ as lattice points in a triangle, for~(a) degree~one, (b)~degree~two, and (c)~degree~three.}\label{fig:triangles}
\end{figure}

Let~$x$, $y$, and~$z$ be the three variables in Question~\ref{quest:intro}. We consider some examples first, and we obtain three monomials of degree one, six monomials of degree two, and ten monomials of degree three; see Figure~\ref{fig:triangles}. In the process, we see that the monomials in~$x$, $y$, and~$z$ of degree~$m$ can be arranged in a triangular shape where the exponent of~$x$ decreases from the top down, and the exponent of~$z$ decreases from the left to the right. The exponent of the remaining variable~$y$ is completely determined by the two other exponents and the total degree~$m$. We conclude that the number of monomials in three variables of a given degree
equals the number of points in a certain triangle, as illustrated by Figure~\ref{fig:triangles}.

More precisely, the number of monomials in~$x,y$, and~$z$ of degree~$m$ equals the number of~\emph{lattice points} -- that is, points with integer-valued coordinates -- in the triangle with vertices~$(0,0)$, $(m,0)$, and~$(m,m)$. In other words, let~$T$ be the triangle with vertices~$(0,0)$, $(1,0)$, and~$(1,1)$, let
\[
mT = \{ m\bv \mid \bv \in T \}
\]
be the~\emph{$m$-th dilation} of~$T$, and let~$L_T(m) = |mT\cap \Z^2|$ count the number of lattice points in~$mT$. Then, the number of monomials in three variables of degree~$m$ equals~$L_T(m)$.

In order to answer Question~\ref{quest:intro}, we thus wish to find a formula for~$L_T(m)$. The set of all lattice points contained in the triangle~$mT$ can be constructed by removing the diagonal from a square of side length~$m+1$ and, subsequently, discarding one of the two resulting congruent triangles.\footnote{Note that a square with side length~$m+1$ contains~$(m+2)^2$ lattice points. The diagonal contains~$m+2$ lattice points.} From this construction, it can be seen that a formula as desired is given by
\[
L_T(m) = \frac{(m+2)^2-(m+2)}{2}=\frac{1}{2}m^2 + \frac{3}{2}m + 1.
\]
Notice that the leading coefficient coincides with the area of~$T$. This is not a coincidence, but part of a bigger story known as~\emph{Ehrhart theory}.

Our goal is to introduce a famous question asked by Tadao~Oda at Oberwolfach in~1997. Roughly speaking, Oda wondered whether the lattice points within a polytope, given that it is of a certain type, always satisfy an elegant counting property; see Question~\ref{Odas-conjecture} at the end of this snapshot for a precise statement. It will become apparent that a certain instance of this problem is related to Question~\ref{quest:intro} in that, in view of how we counted the monomials, it amounts to the question whether each monomial of degree~$m$ can be written as a product of~$m$ variables. Of course, one immediately sees that this is the case. However, Oda's question is not so easy in full generality; although simple to state, it remains unanswered.

Let us now introduce the general picture. We work with the lattice of integral points~$\Z^d \subset \R^d$, that is, the set which consists of all points in~$d$\nobreakdash-dimensional space whose coordinates are integers. A subset~$C\subset\R^d$ is called~\emph{convex} if every straight line segment which connects two points in~$C$ lies entirely within~$C$. The~\emph{convex hull} of a set of points~$B\subset\R^d$ is the inclusion\nobreakdash-wise smallest subset~$\conv(B)\subset\R^d$ which is convex and contains~$B$. For~$B=\{\bv_1,\ldots,\bv_n\}$, an alternative description is as follows:
\[
\conv(\bv_1, \dots, \bv_n) = \left\{ \sum_{i=1}^n \lambda_i \bv_i \ \Bigg|\ \lambda_1, \dots, \lambda_n \in [0,1], \sum_{i=1}^n \lambda_i = 1 \right\}.
\]
The term~\emph{lattice polytope} shall describe the convex hull~$P=\conv(\bv_1,\ldots,\bv_n)$ of a finite number of lattice points~$\bv_1,\ldots,\bv_n$.

As before, the~$m$-th dilation of~$P$ is the polytope~$mP=\{ m \bv \mid \bv \in P\}$. The value of the \emph{Ehrhart function}~$L_P$ at~$m \in \mathbb{Z}_{\ge0}$ is defined to be the number of lattice points~$L_P(m)=|mP \cap \Z^d|$ in~$mP$. The~\emph{dimension} of~$P$, denoted by~$\dim(P)$, is the dimension of the smallest affine space containing~$P$.\footnote{The term~\emph{affine subspace} refers to a point, line, plane,~$\ldots$ in~$\R^d$. For example, a non-degenerate triangle is contained in a plane but in no line; therefore, such a triangle is of dimension two.}

\begin{theorem}[{see~\cite{Ehrhart62, Ehrhart67}}]\label{thm:Ehrhart}
There exists a polynomial~$f$ of degree~$\dim(P)$ with rational coefficients such that~$L_P(m) = f(m)$ for all~$m \in \Z_{\ge0}$. Furthermore, the leading coefficient of~$f$ coincides with the Euclidean volume~$\vol(P)$ of~$P$.
\end{theorem}

Theorem~\ref{thm:Ehrhart} allows us to interpret~$L_P$ as a polynomial of degree~$\dim(P)$ which we call the~\emph{Ehrhart polynomial} of~$P$. That the leading coefficient coincides with the volume of the polytope~$P$ is surprising. However, if~$P$ is a~\emph{polygon}, that is, if~$P$ is a polytope of dimension two, then the relationship between the area of~$P$ and the lattice points it contains can be made even more precise.

\begin{theorem}[Pick's Theorem]
For a lattice polygon~$P \subset \R^2$, the Euclidean volume~$\vol(P)$ is given by the formula
\[
\vol(P) = |P^\circ \cap \Z^2| + \frac{|\partial P \cap \Z^2|}2 - 1.
\]
Here,~$|P^\circ \cap \Z^2|$ denotes the number of interior lattice points, and~$|\partial P \cap \Z^2|$ denotes the number of boundary lattice points of~$P$.
\end{theorem}

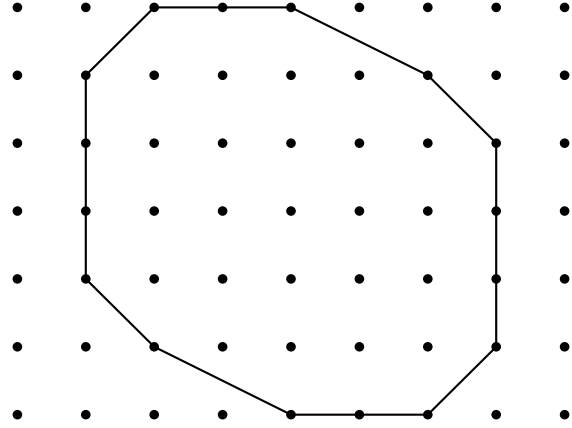
\begin{figure}[t]
\centering
\begin{tikzpicture}[scale=0.9]
\foreach \x in {-1,...,7}
\foreach \y in {0,...,6}
\fill (\x,\y) circle (2pt);
\draw[thick] (3,0) -- (5,0) -- (6,1) -- (6,4) -- (5,5) -- (3,6) -- (1,6) -- (0,5) -- (0,2) -- (1,1) -- cycle;
\end{tikzpicture}
\caption{Computing the area of a lattice polygon by Pick's Theorem. There are~23 interior lattice points and~16 boundary lattice points.\label{fig:ex-pick}}
\end{figure}

Let us use Pick's Theorem to compute the area of the polygon~$P$ shown in Figure~\ref{fig:ex-pick}.
It can be seen that~$P$ has~$23$ interior lattice points and~$16$ boundary lattice points, so its area is given by
\[
\vol(P) = 23 + \frac{16}2 -1 = 30.
\]

\section{Proof of Pick's Theorem}
The proof of Pick's Theorem contains numerous beautiful ideas and constructions from Ehrhart theory. Since several of these play an important role in motivating Oda's Oberwolfach question, we give a proof here. The fundamental idea is to proceed by induction on the number of lattice points~$|P\cap\Z^2|$.\footnote{This means that the theorem is first proved for polygons containing exactly three lattice points (Section~\ref{sec:base-case}). In a second step, it is proved that if the statement of the theorem holds for all polygons which contain at most~$N$ lattice points, then it also holds for all polygons containing~$N+1$ lattice points, where~$N\geq 3$ is any integer (Section~\ref{sec:indictive}). As a consequence, the theorem must then hold for all integers~$N\geq 3$.}

\subsection{The base case}\label{sec:base-case}
We assume that~$P$ is a triangle whose vertices~$\bv_1, \bv_2$, and~$\bv_3$ are its only lattice points; that is,~$|P\cap\Z^2|=3$. Such triangles are called~\emph{empty}. We will show that, in this case,~$\bv_1, \bv_2$, and~$\bv_3$ constitute an~\emph{affine basis} of~$\Z^2$; that is, the difference vectors~$\bv_1 - \bv_3$ and~$\bv_2 - \bv_3$ form a~\emph{lattice basis} of~$\Z^2$.\footnote{The vectors~$\bv_1, \dots, \bv_d \in \Z^d$ are said to form a~\emph{lattice basis} if each vector~$\bv \in \Z^d$ admits a unique representation of the form~$\bv = \lambda_1 \bv_1 + \cdots + \lambda_d \bv_d$ for~$\lambda_1, \dots, \lambda_d \in \Z$; the~$d+1$ elements~$\bv_1, \dots, \bv_{d+1} \in \Z^d$ are said to form an~\emph{affine basis} if the difference vectors~$\bv_1 - \bv_{d+1}, \dots, \bv_d - \bv_{d+1}$ form a lattice basis.} Let us assume this for a moment and show how it implies Pick's Theorem in the base case.

By assumption,~$\bv_1 - \bv_3$ and~$\bv_2 - \bv_3$ form a lattice basis of~$\Z^2$. The $2\times 2$-matrix~$A$ whose rows are these difference vectors is invertible over~$\Z$. In other words, there exists another $2\times 2$-matrix~$B$ with entries in~$\Z$ such that the matrix product~$AB$ equals the identity matrix, and thus the determinant\footnote{The value of the determinant of a $2\times 2$-matrix $\begin{psmallmatrix}a&b\\c&d\end{psmallmatrix}$ is given by the expression~$ad-bc$.} of~$A$ must itself be invertible over~$\Z$; this is, the determinant of~$A$ is~$1$ or~$-1$. It is a geometrical fact that the absolute value of the determinant~$\det\begin{psmallmatrix}\ba\\\bb\end{psmallmatrix}$ of the $2\times 2$-matrix with rows~$\ba,\bb\in\Z^2$ equals the area of the parallelogram spanned by~$\ba$ and~$\bb$. Hence, the area of~$P$ equals~$(1/2)\left|\det(A)\right| = 1/2=0+3/2-1$, as desired.

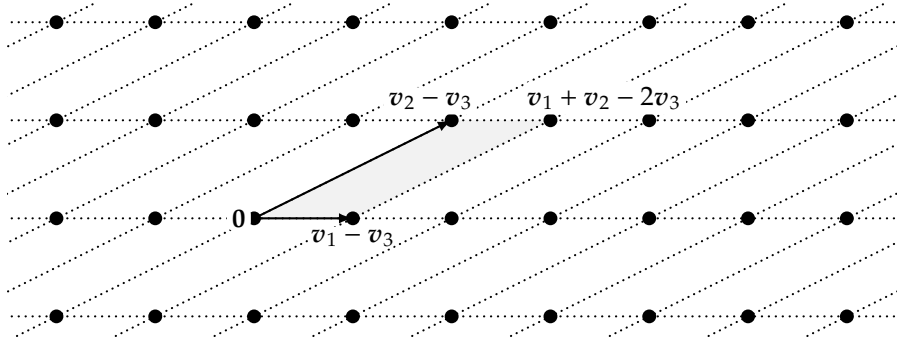
\begin{figure}[ht]
\centering
\begin{tikzpicture}[scale=1.3]
\begin{scope}
\clip (-2.5,-2) rectangle (6.5,2.2);
\fill[fill=gray!10] (0,0) -- (1,0) -- (3,1) -- (2,1) -- cycle;
\draw[thick,dash pattern=on 0pt off 2.5pt, line cap=round] (-4,2) -- (8,2);
\draw[thick,dash pattern=on 0pt off 2.5pt, line cap=round] (-4,1) -- (8,1);
\draw[thick,dash pattern=on 0pt off 2.5pt, line cap=round] (-4,0) -- (8,0);
\draw[thick,dash pattern=on 0pt off 2.5pt, line cap=round] (-4,-1) -- (8,-1);
\foreach \k [%
	evaluate=\k as \xbottom using {-8.4 + \k},%
	evaluate=\k as \xtop using {-1.6 + \k}%
	] in {0,...,14}
		\draw[thick,dash pattern=on 0pt off 2.25pt, line cap=round] (\xbottom,-1.2) -- (\xtop,2.2);

\foreach \x in {-2,...,6}{
	\foreach \y in {-1,...,2}{
		\fill (\x,\y) circle (2pt);
	}
}
\node[left=1pt,fill=white, rounded corners, inner sep=2pt] at (0,0) {$\orig$};
\node[below=3pt,fill=white,rounded corners,inner sep=1pt] at (1,0) {$\bv_1-\bv_3$};
\node[above=2pt,xshift=-8pt,fill=white,rounded corners,inner sep=2pt] at (2,1) {$\bv_2-\bv_3$};
\node[above=2pt,xshift=20pt,fill=white,rounded corners,inner sep=2pt] at (3,1) {$\bv_1+\bv_2-2\bv_3$};
\draw[thick,-latex] (0,0) -- (1,0);
\draw[thick,-latex] (0,0) -- (2,1);
\end{scope}
\end{tikzpicture}
\vspace{-2em}
\caption{Tiling the plane with translations of a parallelogram which comes from an empty triangle with vertices~$\bv_1,\bv_2,\bv_3$.\label{fig:tiling}}
\end{figure}

It remains to show that~$\bv_1, \bv_2$, and~$\bv_3$ together form an affine basis of~$\Z^2$. Let~$\cP$ denote the parallelogram spanned by~$\bv_1 - \bv_3$ and~$\bv_2 - \bv_3$; that is,
\[
\cP = \{ \lambda_1 \bv_1 + \lambda_2 \bv_2 - (\lambda_1 + \lambda_2)\bv_3 \mid 0 \le \lambda_i \le 1 \}.
\]
Such a parallelogram is shown in Figure~\ref{fig:tiling}.
We claim that the only lattice points contained in~$\cP$ are its vertices~$\orig,\bv_1-\bv_3,\bv_2-\bv_3$, and~$\bv_1+\bv_2-2\bv_3$. Indeed, the triangle with the vertices~$\orig,\bv_1-\bv_3,$ and~$\bv_2-\bv_3$ is empty, and so its opposite triangle, that is, the triangle with the vertices~$\bv_1-\bv_3,\bv_2-\bv_3$, and~$\bv_1+\bv_2-2\bv_3$ which covers the other half of~$\cP$, is also empty. Certainly, the set~$\{\cP+\bv\mid\bv\in\Z^2\}$ of all translations of~$\cP$ tiles the plane, as illustrated in Figure~\ref{fig:tiling}, from which we may observe that every lattice point in~$\Z^2$ is a vertex of a lattice translation of~$\cP$. Hence, each~$\bx\in\Z^2$ can be expressed as a linear combination~$\bx=k_1(\bv_1-\bv_3)+k_2(\bv_2-\bv_3)$ for~$k_1,k_2\in\Z$. Moreover, this representation is unique, so~$\bv_1,\bv_2,$ and~$\bv_3$ together form an affine basis of~$\Z^2$.

\subsection{The inductive step}\label{sec:indictive}
Let us assume Pick's Theorem is valid for all lattice polygons~$Q \subset \R^2$ with~$|Q\cap \Z^2|\le N$ for some~$N\ge3$. We need to show that it then holds for all lattice polygons~$P$ with~$|P\cap \Z^2|=N+1$. Let~$P$ be such a polygon. The inductive step splits into two cases: either there are exactly three lattice points on the boundary~$\partial P$ of~$P$; or~$|\partial P\cap\Z^2|>3$.

\begin{figure}[t]
	\centering
	\begin{subfigure}{.3\textwidth}
		\centering
		\begin{tikzpicture}[scale=1]
			\foreach \x in {-1,...,2}
			\foreach \y in {-1,...,3}
			\fill (\x,\y) circle (2pt);
			\draw[thick] (2,0) -- (0,3) -- (-1,-1) -- cycle;
			\draw[thin,dashed] (0,0) -- (2,0);
			\draw[thin,dashed] (0,0) -- (0,3);
			\draw[thin,dashed] (0,0) -- (-1,-1);
			\node[above right] at (0,0) {$\bv$};
		\end{tikzpicture}
		\vspace{30pt}
		\caption{\label{subfig:split-in-three-left}}
	\end{subfigure}
	\begin{subfigure}{.5\textwidth}
		\centering
		\begin{tikzpicture}[scale=1]
			\draw[very thick,decorate,decoration={zigzag,amplitude=.5mm, segment length=7mm, post length=0mm},draw=lightgray!75] (0.75,.125) -- (4,.125);
			\draw[very thick,decorate,decoration={zigzag,amplitude=.5mm, segment length=7mm, post length=0mm},draw=lightgray!75] (0.5,.5) -- (0.5,4);
			\draw[very thick,decorate,decoration={zigzag,amplitude=.5mm, segment length=7mm, post length=0mm},draw=lightgray!75] (0.25,0) -- (-1.5,-1.75);

			\begin{scope}[shift={(1,.5)}]
				\draw[thick] (2,0) -- (0,3);
				\draw[thin,dashed] (0,3) -- (0,0) -- (2,0);
				\node (P1) at (3,1.5) {$P_1$};
				\draw[very thin] (P1.west) edge[-latex,bend right] (1,1.5);
				\foreach \x in {0,...,2}
				\foreach \y in {0,...,3}
				\fill[fill=black] (\x,\y) circle (2pt);
				\fill[fill=black] (0,0) circle (2pt);
			\end{scope}

			\begin{scope}[shift={(0,.25)}]
				\draw[thick] (0,3) -- (-1,-1);
				\draw[thin,dashed] (0,3) -- (0,0) -- (-1,-1);
				\node (P2) at (-2,1.5) {$P_2$};
				\draw[very thin] (P2.east) edge[-latex,bend left] (-.5,1);
				\foreach \x in {-1,0}
				\foreach \y in {0,...,3}
				\fill[fill=black] (\x,\y) circle (2pt);
				\fill[fill=black] (-1,-1) circle (2pt);
			\end{scope}

			\begin{scope}[shift={(0.5,-.25)}]
				\draw[thick] (2,0) -- (-1,-1);
				\draw[thin,dashed] (2,0) -- (0,0) -- (-1,-1);
				\node (P3) at (2,-1.5) {$P_3$};
				\draw[very thin] (P3.west) edge[-latex,bend left] (.5,-.5);
				\foreach \x in {0,...,2}
				\foreach \y in {-1,0}
				\fill (\x,\y) circle (2pt);
				\fill (-1,-1) circle (2pt);
			\end{scope}
	\end{tikzpicture}
	\caption{\label{subfig:split-in-three-right}}
	\end{subfigure}
	\caption{A lattice polygon~\ref*{subfig:split-in-three-left} with exactly three boundary lattice points is cut~\ref*{subfig:split-in-three-right} into three subpolygons~$P_1,P_2$, and~$P_3$.
	Here,~$b_1=4$, $b_2=3$, $b_3=2$ and~$i_1=1$, $i_2=i_3=0$.\label{fig:ind_step_boundary_three_lattice_pts}}
\end{figure}
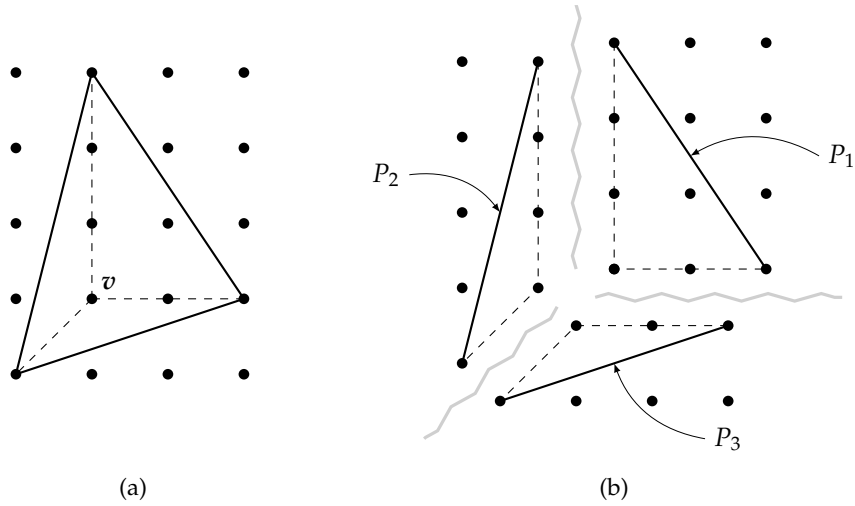

If~$|\partial P\cap\Z^2|=3$ then the interior~$P^\circ$ of~$P$ contains at least one lattice point~$\bv$. The three line segments connecting~$\bv$ with the three vertices of~$P$ split the polygon into three subpolygons~$P_1$, $P_2$, and~$P_3$, as illustrated by Figure~\ref{fig:ind_step_boundary_three_lattice_pts}. We shall count the lattice points in each of the~$P_k$, $k=1,2,3$, and relate these counts to the number of lattice points in~$P$. Let~$i_k=|P_k^\circ\cap\Z^2|$, and let~$b_k=|\partial P_k\cap P^\circ\cap\Z^2|=|\partial P_k\cap\Z^2|-2$, for~$k=1,2,3$.

Since~$|P_k\cap\Z^2|\leq N$ for~$k=1,2,3$, we can apply our initial assumption on the validity of Pick's Theorem for lattice polygons containing at most~$N$ lattice points to each of the polygons~$P_1,P_2,P_3$, and we obtain
\begin{align*}
\vol(P_k) &= |P_k^\circ \cap \Z^2| + \frac{|\partial P_k \cap \Z^2|}2 -1\\
&= i_k + \frac{b_k+2}2-1.
\end{align*}
The area of~$P$ equals the sum of the areas of the subpolygons~$P_1,P_2$, and~$P_3$. By further counting the lattice points in each of the three subpolygons separately, we arrive at the following chain of equalities:
\begin{align*}
\vol(P) &= \sum_{k=1}^3 \vol(P_k)\\
&= \sum_{k=1}^3 \left(i_k+\frac{b_k+2}2-1\right)\\
&= \sum_{k=1}^3 i_k + \frac{\sum_{k=1}^3 b_k}2\\
&= |P^\circ\cap\Z^2| + \frac12\\
&= |P^\circ\cap\Z^2| + \frac32 - 1\\
&= |P^\circ \cap \Z^2| + \frac{|\partial P \cap \Z^2|}2 - 1.
\end{align*}
Notice that each lattice point which lies on the interior parts of the dissecting lines is counted twice by the term~$\sum_{k=1}^3b_k$, except for~$\bv$, which is counted three times. This is accounted for by the addition of~$1/2$ in the fourth line.

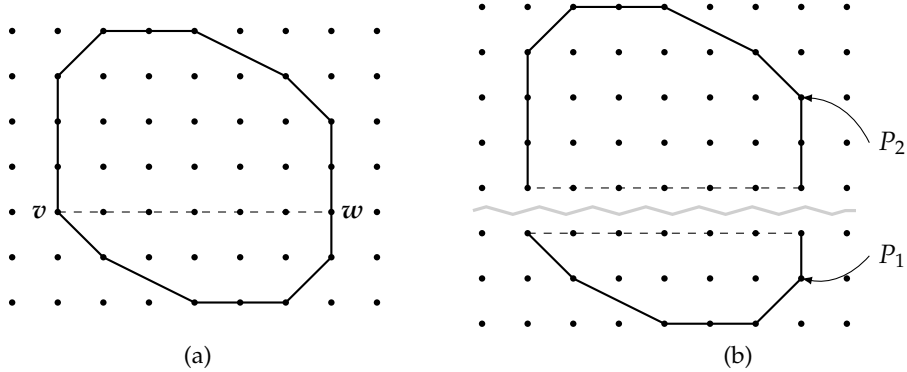
\begin{figure}[!ht]
	\begin{subfigure}{.45\textwidth}
		\centering
		\begin{tikzpicture}[scale=0.6]
			\begin{scope}
				\foreach \x in {-1,...,7}
				\foreach \y in {0,...,6}
				\fill (\x,\y) circle (2pt);
				\draw[thick] (3,0) -- (5,0) -- (6,1) -- (6,4) -- (5,5) -- (3,6) -- (1,6) -- (0,5) -- (0,2) -- (1,1) -- cycle;
				\draw[dashed] (0,2) node[left] {$\bv$} -- (6,2) node[right] {$\bw$};
			\end{scope}
		\end{tikzpicture}
		\vspace{8pt}
		\caption{\label{subfig:split-in-two-left}}
	\end{subfigure}
	\begin{subfigure}{.45\textwidth}
		\begin{tikzpicture}[scale=0.6]
			\draw[very thick,decorate,decoration={zigzag,amplitude=.5mm, segment length=7mm, post length=0mm},draw=lightgray!75] (-1.2,2) -- (7.2,2);

			\begin{scope}[shift={(0,0.5)}]
				\foreach \x in {-1,...,7}
				\foreach \y in {2,...,6}
				\fill (\x,\y) circle (2pt);
				\draw[thick] (6,2) -- (6,4) -- (5,5) -- (3,6) -- (1,6) -- (0,5) -- (0,2);
				\draw[dashed] (6,2) -- (0,2);
			\end{scope}

			\begin{scope}[shift={(0,-0.5)}]
				\foreach \x in {-1,...,7}
				\foreach \y in {0,...,2}
				\fill (\x,\y) circle (2pt);
				\draw[thick] (0,2) -- (1,1) -- (3,0) -- (5,0) -- (6,1) -- (6,2);
				\draw[dashed] (0,2) -- (6,2);

				\node (P2) at (8,4) {$P_2$};
				\draw (P2.west) edge[-latex,bend right] (6,5);
				\node (P1) at (8,1.5) {$P_1$};
				\draw (P1.west) edge[-latex,bend left] (6,1);
			\end{scope}
		\end{tikzpicture}
		\caption{\label{subfig:split-in-two-right}}
	\end{subfigure}
	\caption{A lattice polygon~\ref*{subfig:split-in-two-left} with four or more lattice points on its boundary can be split~\ref*{subfig:split-in-two-right} into two subpolygons~$P_1$ and~$P_2$.
	Here,~$i_1=4$, $i_2=14$,~$b_1=7$, $b_2=11$, and~$i=5$.\label{fig:pick_ind_step_boundary_good}}
\end{figure}

If~$|\partial P\cap\Z^2|>3$, we can find~$\bv, \bw \in \partial P \cap \Z^2$ such that the line segment connecting~$\bv$ and~$\bw$ divides~$P$ into two subpolygons~$P_1$ and~$P_2$ with~$|P_k\cap\Z^2|\leq N$ for~$k=1,2$. This is illustrated by Figure~\ref{fig:pick_ind_step_boundary_good}.
We again count the number of lattice points in the subpolygons~$P_1$ and~$P_2$ and relate these counts to the number of lattice points in~$P$. Let~$i_k=|P_k^\circ\cap\Z^2|$, and let~$b_k=|(\partial P_k\setminus P^\circ)\cap\Z^2|$ be the number of boundary lattice points of~$P_k$ that are not contained in the interior of~$P$. By~$i$ we denote the number of interior lattice points in the line segment from~$\bv$ to~$\bw$. Since~$|P_k\cap\Z^2|\leq N$, we can apply the inductive hypothesis to obtain
\[
\vol(P_k) = i_k + \frac{b_k+i}2 - 1.
\]
The area of~$P$ equals the sum of the areas of~$P_1$ and~$P_2$; hence,
\begin{align*}
\vol(P) &= \vol(P_1)+\vol(P_2)\\
&=\left(i_1+\frac{b_1+i}{2}-1\right)+\left(i_2+\frac{b_2+i}{2}-1\right)\\
&= (i_1 + i_2+i) + \frac{b_1+b_2}{2} -2\\
&= |P^\circ \cap \Z^2| + \frac{|\partial P \cap \Z^2|}2 -1.
\end{align*}
Regarding the final equality, notice that~$i_1+i_2+i=|P^\circ\cap\Z^2|$, whereas~$b_1+b_2$ counts the lattice points on the boundary of~$P$, counting~$\bv$ and~$\bw$ twice, however.

\section{Oda's Oberwolfach question}
Recall our reinterpretation of Question~\ref{quest:intro} as a problem on the number of lattice points in dilations of the empty triangle~$T$ with vertices at~$(0,0),(1,0)$, and~$(1,1)$; the lattice points contained in the set~$T\cap\Z^2$ each correspond to a monomial of degree one. We now construct a single combinatorial object which simultaneously encodes the lattice points (and hence their counts) in all dilations~$mT$,~$m \in \Z_{\ge0}$.

\begin{figure}[!ht]
\centering

\makeatletter
\pgfdeclareradialshading[tikz@ball]{ball}{\pgfqpoint{-20bp}{20bp}}{%
color(0bp)=(tikz@ball!30!white); 
color(13bp)=(tikz@ball!0!gray); 
color(21bp)=(tikz@ball!70!black); 
color(25bp)=(black!70); 
color(30bp)=(black!70)} 
\makeatother 

\begin{tikzpicture}[%
rotate around y=115,%
rotate around z=0,%
scale=.7,%
line cap=round,%
line join=round%
]
\fill[fill=gray,opacity=.35] (0,0,0) -- (4.5,4.5,4.5) -- (0,4.5,4.5) -- cycle;
\fill[color=gray,opacity=.15] (0,0,0) -- (0,4.5,0) -- (4.5,4.5,4.5);
\draw[thick,-latex] (0,0,-1.5) -- (0,0,4.5) node[below] {$x$};
\draw[thick,-latex] (-1.5,0,0) -- (4.5,0,0) node[left] {$y$};
\draw[thick,-latex] (0,-1,0) -- (0,5.5,0) node[left] {$z$};
\node[right=.25cm] at (0,1,1) {$T\times\{1\}$};
\node[right=.25cm] at (0,2,2) {$2T\times\{2\}$};
\node[right=.25cm] at (0,3,3) {$3T\times\{3\}$};
\fill[fill=gray,opacity=.5] (0,1,0) -- (0,1,1) -- (1,1,1) -- cycle;
\draw[thick,color=gray!90] (0,1,0) -- (0,1,1);
\draw[dash pattern=on 2pt off 3pt,thick,color=gray!120] (0,1,0) -- (1,1,1) -- (0,1,1);
\fill[fill=gray,opacity=.5] (0,2,0) -- (0,2,2) -- (2,2,2) -- cycle;
\draw[thick,color=gray] (0,2,0) -- (0,2,2);
\draw[dash pattern=on 2pt off 3pt,thick,color=gray] (0,2,0) -- (2,2,2) -- (0,2,2);
\fill[fill=gray,opacity=.5] (0,3,0) -- (0,3,3) -- (3,3,3) -- cycle;
\draw[thick,color=gray] (0,3,0) -- (0,3,3);
\draw[dash pattern=on 2pt off 3pt,thick,color=gray] (0,3,0) -- (3,3,3) -- (0,3,3);
\draw[thick,color=gray] (0,0,0) -- (0,5,0);
\draw[dashed,thick,color=gray] (0,0,0) -- (4.5,4.5,4.5);
\draw[thick,color=gray] (0,0,0) -- (0,5,5);
\draw[thick,color=gray] (4.5,4.5,4.5) -- (5,5,5);
\fill[fill=gray,opacity=.5] (0,0,0) -- (0,4.5,0) -- (0,4.5,4.5) -- cycle;
\shade[ball color=darkgray] (0,1,0) circle (2pt);
\shade[ball color=darkgray] (0,1,1) circle (2pt);
\shade[ball color=darkgray] (1,1,1) circle (2pt);
\shade[ball color=darkgray] (0,2,0) circle (2pt);
\shade[ball color=darkgray] (0,2,1) circle (2pt);
\shade[ball color=darkgray] (1,2,1) circle (2pt);
\shade[ball color=darkgray] (1,2,2) circle (2pt);
\shade[ball color=darkgray] (0,2,2) circle (2pt);
\shade[ball color=darkgray] (2,2,2) circle (2pt);
\shade[ball color=darkgray] (0,3,0) circle (2pt);
\shade[ball color=darkgray] (0,3,1) circle (2pt);
\shade[ball color=darkgray] (1,3,1) circle (2pt);
\shade[ball color=darkgray] (1,3,2) circle (2pt);
\shade[ball color=darkgray] (0,3,2) circle (2pt);
\shade[ball color=darkgray] (2,3,2) circle (2pt);
\shade[ball color=darkgray] (0,3,3) circle (2pt);
\shade[ball color=darkgray] (3,3,3) circle (2pt);
\shade[ball color=darkgray] (1,3,3) circle (2pt);
\shade[ball color=darkgray] (2,3,3) circle (2pt);
\end{tikzpicture}
\caption{The cone~$C_T$ over the empty triangle~$T$ affinely embedded into~$\R^3$.\label{fig:cone_over_empty_triangle}}
\end{figure}
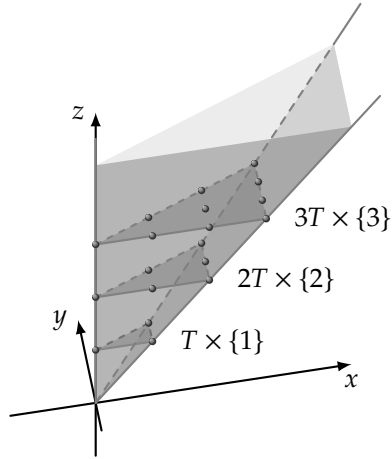

We begin by identifying~$\R^2$ with the subset~$\R^2\times\{1\}$ of~$\R^3$; in other words, we think of the plane as lying inside three\nobreakdash-dimensional space at height one. In this copy of~$\R^2$, we consider the triangle~$T$. The~\emph{cone over~$T$} is defined as~$C_T = \{ \lambda (\bx, 1) \mid \lambda \in \R_{\ge0}, \bx \in T\}$; it is depicted in Figure~\ref{fig:cone_over_empty_triangle}. We may think of~$C_T$ as the union of all half-lines that start at the origin and pass through a point in~$T\times\{1\}$. The cross-section of~$C_T$ at height~$m\in\Z_{\ge0}$, or, more precisely,~$C_T\cap\{(x,y,z)\mid z=m\}=(mT,m)$, can be identified with the dilation~$mT$ of~$T$. Hence, the cone~$C_T$ encodes the lattice points in all dilations of~$T$.

This construction can be generalised to arbitrary lattice polytopes~$P \subset \R^d$. Indeed, the~\emph{cone over~$P$} can be defined as
\begin{align*}
C_P &= \R_{\ge0} (P\times\{1\})\\
&= \{ \lambda (\bx, 1) \mid \lambda \in \R_{\ge0}, \bx \in P\}.
\end{align*}
As before, the cross-section of~$C_P$ at height~$m\in\Z_{\ge0}$ can be identified with the dilation~$mP$.

\subsection{The integer decomposition property}
A lattice polytope~$P\subset\R^d$ is said to have the~\emph{integer decomposition property} if each lattice point in~$C_P$ at height~$m$ can be written as the sum of~$m$ (not necessarily distinct) lattice points at height one.

\begin{example}
If the polytope under consideration is an empty triangle~$T$, this means that for all $(\bx,m)\in C_T\cap\Z^3$ there exist (not necessarily distinct) points~$(\bx_1,1),\ldots,(\bx_m,1)\in C_T\cap\Z^3$ such that
\[
(\bx,m)=(\bx_1,1)+\cdots+(\bx_m,1).
\]
Indeed, every empty triangle~$T$ has the integer decomposition property; the interested reader is invited to think about why this is the case.
\end{example}

It is natural to ask which other lattice polytopes have the integer decomposition property. Let us study some important classes of lattice polytopes for which this question has been resolved.

\begin{example}
We generalise empty lattice triangles to arbitrary dimension.
If~$\bv_0, \dots, \bv_d$ form an affine basis of~$\Z^d$, their convex hull~$S = \conv(\bv_0, \dots, \bv_d)$ is called a~\emph{$d$\nobreakdash-dimensional unimodular simplex}. We claim that such a simplex~$S$ has the integer decomposition property.
Indeed, since~$\bv_0, \dots, \bv_d$ is an affine basis of~$\Z^d$, it follows that~$(\bv_0,1), \dots, (\bv_d,1)$ is a lattice basis of~$\Z^{d+1}$. Hence, every lattice point~$(\bx,m)$ in the cone~$C_S$ can be expressed as a linear combination~$(\bx,m) = \sum_{i=0}^d\lambda_i (\bv_i,1)$, where~$\lambda_i \in \Z$.
Since~$(\bx,m)$ is in~$C_S$, it follows that~$\lambda_i \ge 0$, and thus~$S$ has the integer decomposition property.
\end{example}

\begin{figure}[t]
\centering
\begin{tikzpicture}[scale=0.8]
\foreach \x in {-1,...,7}
\foreach \y in {0,...,6}
\fill (\x,\y) circle (2pt);
\draw[thick] (3,0) -- (5,0) -- (6,1) -- (6,4) -- (5,5) -- (3,6) -- (1,6) -- (0,5) -- (0,2) -- (1,1) -- cycle;
\draw[very thin] (0,5) -- (5,5);
\draw[very thin] (0,4) -- (6,4);
\draw[very thin] (0,3) -- (6,3);
\draw[very thin] (0,2) -- (6,2);
\draw[very thin] (1,1) -- (6,1);
\draw[very thin] (1,6) -- (1,1);
\draw[very thin] (2,6) -- (2,1);
\draw[very thin] (3,6) -- (3,0);
\draw[very thin] (4,5) -- (4,0);
\draw[very thin] (5,5) -- (5,0);
\draw[very thin] (0,4) -- (2,6);
\draw[very thin] (0,3) -- (3,6);
\draw[very thin] (0,2) -- (3,5);
\draw[very thin] (3,6) -- (4,5);
\draw[very thin] (1,2) -- (4,5);
\draw[very thin] (1,1) -- (5,5);
\draw[very thin] (2,1) -- (3,0);
\draw[very thin] (2,1) -- (5,4);
\draw[very thin] (3,1) -- (6,4);
\draw[very thin] (3,0) -- (6,3);
\draw[very thin] (4,0) -- (6,2);
\draw[very thin] (5,0) -- (6,1);
\end{tikzpicture}
\caption{Every two\nobreakdash-dimensional polygon can be covered by empty triangles.\label{fig:triangulation}}
\end{figure}
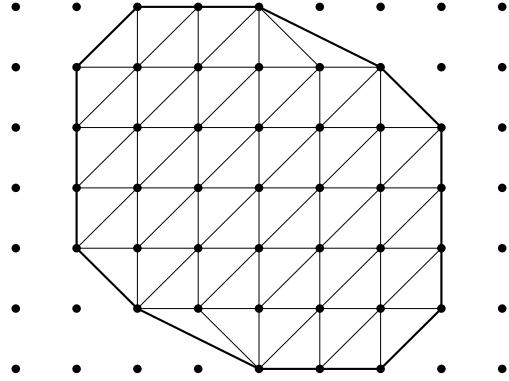

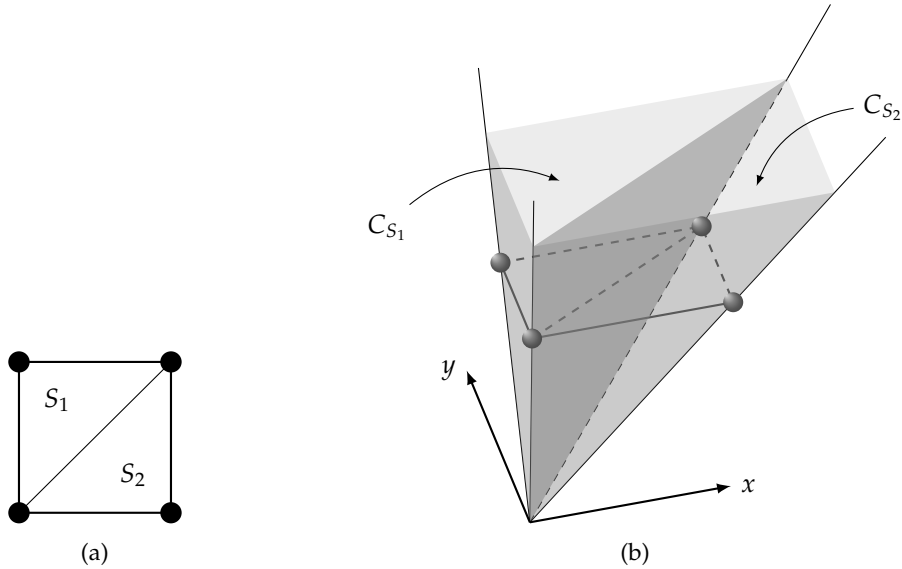
\begin{figure}[t]
	\centering
	\begin{subfigure}{.3\textwidth}
		\centering
		\begin{tikzpicture}[scale=2]
			\begin{scope}
				\foreach \x in {0,1} \foreach \y in {0,1}
				\fill[fill=black] (\x,\y) circle (2pt);
				\draw[thick] (0,0) -- (1,0) -- (1,1) -- (0,1) -- cycle;
				\draw[thin] (0,0) -- (1,1);
				\node at (.25,.75) {$S_1$};
				\node at (.75,.25) {$S_2$};
			\end{scope}
		\end{tikzpicture}
		\caption{\label{subfig:triangulation-polytope}}
	\end{subfigure}
	\begin{subfigure}{.6\textwidth}
		\centering
		\makeatletter
		\pgfdeclareradialshading[tikz@ball]{ball}{\pgfqpoint{-20bp}{20bp}}{%
		color(0bp)=(tikz@ball!30!white); 
		color(13bp)=(tikz@ball!0!gray); 
		color(21bp)=(tikz@ball!70!black); 
		color(25bp)=(black!70); 
		color(30bp)=(black!70)} 
		\makeatother 
		\begin{tikzpicture}
			\begin{scope}[
			shift={(0,-1)},%
			rotate around y=120,%
			rotate around z=4,%
			scale=2.5%
			]
				\coordinate (O) at (0,0,0);

				\coordinate (A) at (0,1,0);
				\coordinate (B) at (1,1,0);
				\coordinate (C) at (1,1,1);
				\coordinate (D) at (0,1,1);

				\coordinate (A high) at (0,1.5,0);
				\coordinate (B high) at (1.5,1.5,0);
				\coordinate (C high) at (1.5,1.5,1.5);
				\coordinate (D high) at (0,1.5,1.5);

				\coordinate (A top) at (0,1.75,0);
				\coordinate (B top) at (1.75,1.75,0);
				\coordinate (C top) at (1.75,1.75,1.75);
				\coordinate (D top) at (0,1.75,1.75);

				\draw[thick,color=darkgray] (D) -- (A) -- (B);
				\draw[thick,dashed,color=darkgray] (B) -- (C) -- (D);
				\draw (O) -- (A top);
				\draw (O) -- (B top);
				\draw[dashed] (O) -- (C high);
				\draw (C high) -- (C top);
				\draw (O) -- (D top);

				\draw[thick,dashed,color=darkgray] (A) -- (C);

				\fill[color=gray,opacity=.3] (O) -- (A high) -- (B high) -- cycle;
				\fill[color=gray,opacity=.15] (O) -- (B high) -- (C high) -- cycle;
				\fill[color=gray,opacity=.15] (O) -- (C high) -- (D high) -- cycle;
				\fill[color=gray,opacity=.3] (O) -- (D high) -- (A high) -- cycle;
				\fill[color=gray,opacity=.5] (O) -- (A high) -- (C high) -- cycle;
				\draw[-latex,thick] (O) -- (2,0,0) node[left] {$y$};
				\draw[-latex,thick] (O) -- (0,0,1) node[right] {$x$};
				\shade[ball color=gray] (A) circle (1.5pt);
				\shade[ball color=gray] (B) circle (1.5pt);
				\shade[ball color=gray] (C) circle (1.5pt);
				\shade[ball color=gray] (D) circle (1.5pt);

				\draw (.5,1.5,-.5) node[left] {$C_{S_1}$} edge[-latex,bend left] (.75,1.5,.25);
				\draw (1,1.5,1.75) node[right] {$C_{S_2}$} edge[-latex,bend right] (.25,1.5,1.15);
			\end{scope}
		\end{tikzpicture}
		\caption{\label{subfig:triangulation-cone}}
	\end{subfigure}
	\caption{If~\ref*{subfig:triangulation-polytope} triangles~$S_i$ cover a polygon~$P$, then~\ref*{subfig:triangulation-cone} the cones~$C_{S_i}$ cover~$C_P$.\label{fig:triangulation_covers_cone}}
\end{figure}

\begin{example}
Every two\nobreakdash-dimensional lattice polygon can be covered by empty lattice triangles; Figure~\ref{fig:triangulation} provides an example of how this can be achieved. More generally, suppose~$P$ is a lattice polytope in~$\R^d$ which is covered by unimodular simplices~$S_1, \dots, S_n$; that is,~$P = S_1 \cup \cdots \cup S_n$. Such a covering is called a~\emph{unimodular covering} of~$P$. The cone~$C_P$ over~$P$ is covered by the cones~$C_{S_1}, \dots, C_{S_n}$ over the simplices~$S_i$; Figure~\ref{fig:triangulation_covers_cone} illustrates this for the case~$d=2$. Hence, every lattice point in~$C_P$ at height~$m$ can be written as a sum of~$m$ lattice points at height one, namely those which are located at height one in the corresponding cone~$C_{S_i}$. In particular, all two\nobreakdash-dimensional lattice polygons have the integer decomposition property, since they admit a unimodular covering. In higher dimensions it is typically very hard to determine whether a unimodular covering even exists.\footnote{The unimodular coverings in Figures~\ref{fig:triangulation} and~\ref{fig:triangulation_covers_cone} have the property that the simplices intersect along common faces. Such coverings are called~\emph{triangulations}. In dimension two, a unimodular triangulation always exists. However, there exist higher-dimensional lattice polytopes that admit a unimodular covering but no unimodular triangulation~\cite[Example~10]{FirlaZiegler99}.}
\end{example}

\begin{example}
Reeve~\cite{Ree57} describes an infinite family of tetrahedra whose members we now call the~\emph{Reeve tetrahedra}; see Figure~\ref{fig:reeve}. Let
\[
R_r=\conv((0, 0, 0), (1, 0, 0), (0, 1, 0), (1, 1, r))
\]
for some~$r\in\Z_{\geq 1}$.
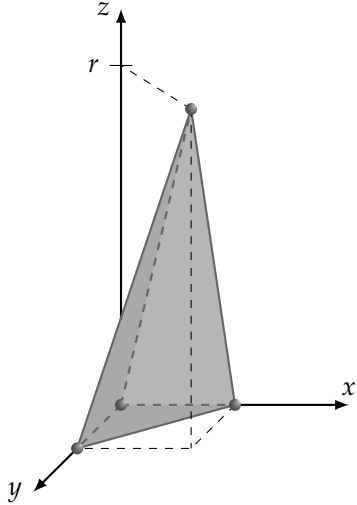
\begin{figure}[t]
\centering
\makeatletter
\pgfdeclareradialshading[tikz@ball]{ball}{\pgfqpoint{20bp}{20bp}}{%
color(0bp)=(tikz@ball!30!white); 
color(13bp)=(tikz@ball!0!gray); 
color(21bp)=(tikz@ball!70!black); 
color(25bp)=(black!70); 
color(30bp)=(black!70)} 
\makeatother
\begin{tikzpicture}[scale=1.5]
\draw[thick,-latex] (0,.75,0) -- (0,3.5,0) node[left] {$z$};
\draw[thick,-latex] (1,0,0) -- (2,0,0) node[above] {$x$};
\draw[thick,-latex] (0,0,1) -- (0,0,2) node[left] {$y$};
\fill[fill=gray,opacity=.15] (0,0,0) -- (1,0,0) -- (1,3,1) -- cycle;
\fill[fill=gray,opacity=.35] (0,0,0) -- (0,0,1) -- (1,3,1) -- cycle;
\fill[fill=gray,opacity=.35] (0,0,0) -- (1,0,0) -- (0,0,1) -- cycle;
\fill[fill=gray,opacity=.5] (1,0,0) -- (1,3,1) -- (0,0,1) -- cycle;

\draw[dashed,thick,draw=gray!120] (0,0,0) -- (1,3,1);
\draw[dashed,thick,draw=gray!120] (0,0,0) -- (1,0,0);
\draw[dashed,thick,draw=gray!120] (0,0,0) -- (0,0,1);
\draw[thick,draw=gray!120] (1,0,0) -- (1,3,1);
\draw[thick,draw=gray!120] (0,0,1) -- (1,3,1);
\draw[thick,draw=gray!120] (1,0,0) -- (0,0,1);

\draw[thin,dashed] (1,3,1) -- (1,0,1);
\draw[thin,dashed] (1,0,0) -- (1,0,1) -- (0,0,1);

\draw[thin,dashed] (1,3,1) -- (0,3,0);

\shade[ball color=gray] (0,0,0) circle (1.5pt);
\shade[ball color=gray] (1,0,0) circle (1.5pt);
\shade[ball color=gray] (0,0,1) circle (1.5pt);
\shade[ball color=gray] (1,3,1) circle (1.5pt);

\draw (-.1,3,0) node[left] {$r$} -- (.1,3,0);
\end{tikzpicture}
\caption{The Reeve tetrahedron~$R_r$.\label{fig:reeve}}
\end{figure}
The Reeve tetrahedron~$R_r$ contains exactly four lattice points, has no interior lattice points, and~$\vol(R_r)=r/6$. Moreover, the number of lattice points in its~$m$-th dilation is given by
\[
L_{R_r}(m)=\frac{r}{6}m^3+m^2+\left(2-\frac{r}{6}\right)m+1.
\]

It is readily seen that~$R_r$ does not admit a unimodular covering when~$r>1$, as its four vertices do not form an affine lattice basis in this case. This stands in contrast with the two\nobreakdash-dimensional case, where every empty triangle is unimodular. The existence of empty simplices whose vertices do not form an affine basis is what makes the study of higher\nobreakdash-dimensional lattice polytopes much richer -- and harder -- than that of polygons. In particular, we should not expect a direct analogue of Pick's Theorem in dimension three or more.

We now show that~$R_r$ does not satisfy the integer decomposition property. Let~$C_r$ be the cone over~$R_r$. The set of lattice points~$C_r\cap\Z^4$ forms a semigroup with~$\Z_{\geq 0}$-basis given by
\[
\begin{cases}
\parbox{6.2cm}{$ (0,0,0,1), (1,0,0,1), (0,1,0,1), (1,1,1,1)$}&\text{ if }r=1,\text{ or}\\[1ex]
\parbox{6.2cm}{\raggedleft$ (0,0,0,1), (1,0,0,1), (0,1,0,1), (1,1,r,1),$\\$(1,1,1,2), \dots, (1,1,r-1,2) \hphantom{,}$} & \text{ if }r>1.
\end{cases}
\]
In particular, for~$r>1$ the point~$(1,1,1,2)\in C_r$ at height two cannot be written as the sum of two points at height one.
\end{example}

\subsection{The smooth case}
Let~$\be_i\in\R^d$ denote the~$i$-th standard basis vector, and let us write
\[
\cone(X) = \left\{ \sum_{i=1}^n \lambda_i \bx_i \ \Bigg|\ n \in \Z_{\ge0}, \lambda_i \in \R_{\ge0}, \bx_i \in X\right\}
\]
for the convex cone generated by a (not necessarily finite) subset~$X\subset\R^d$. A cone~$C\subset\R^d$ is called~\emph{smooth} if it can be identified with~$\cone(\be_1, \dots, \be_d)$ via a change of basis of~$\Z^d$. In other words, it is called smooth if there exists a lattice basis~$\bb_1,\ldots,\bb_d$ of~$\Z^d$ with the property that the map~$\Z^d\rightarrow\Z^d$ sending~$(x_1,\ldots,x_d)\mapsto\sum_{k=1}^dx_k\bb_k\in\Z^d$ is invertible and sends~$\cone(\be_1,\ldots,\be_d)$ to~$C$.

A lattice polytope~$P \subset \R^d$ is called~\emph{smooth} if~$\cone(P - \bv)$ is smooth for every vertex~$\bv$ of~$P$. Oda's Oberwolfach question now asks whether smoothness is a sufficient condition for a lattice polytope to have the integer decomposition property.

\begin{question}[{Oda's Oberwolfach Question~\cite{Oda}}]\label{Odas-conjecture}
Does every smooth lattice polytope have the integer decomposition property?
\end{question}

An indication that this may indeed be the case lies in the fact that the smoothness property of a polytope~$P$ ensures that each of its corner is covered by a dilation of a unimodular simplex, and in the speculation that it should be possible to extend these corner covers sufficiently far inside so as to yield a unimodular covering of the polytope~$P$.

Despite many efforts to answer Question~\ref{Odas-conjecture}, it is still wide open. Moreover, there is currently no general consensus on the likely answer. Substantial efforts have been made in order to find a counterexample; for instance, such an effort is due to Bruns~\cite{Bruns13}. Meanwhile, the study of Oda's Oberwolfach question has led to a considerable number of beautiful results on lattice polytopes which were obtained by a variety of authors~\cite{BeckEtAl19,Gubeladze,HaaseHofmann, HaaseEtAl21}.
\bibliographystyle{plain}
\bibliography{bibliography}
\end{document}